\documentclass[11pt]{article}
\usepackage{latexsym,amsfonts,euscript}

\def\la{\langle}\def\ra{\rangle}

\def\pf{\noindent{\bf Proof\quad }}
\def\pfend{\hfill{$\Box$}}

\newtheorem{lem}{ \bf Lemma}[section]

\newtheorem{pro}[lem]{\bf Proposition}

\title{A CHARACTERIZATION  OF $L_2(2^f)$ IN TERMS OF CHARACTER ZEROS
\thanks{Project supported by  the NNSF of China (Grant No. 10171074), and the NSF of Jiangsu Provincial Education
Department(03KJB110002). }}

\author{Guohua Qian\\
{\footnotesize\small Dept. Math., Changshu Institute of
Technology,
Changshu 215500, Jiangsu, P. R. China.}\\
\footnotesize\small E-mail: ghqian2000@yahoo.com.cn\\
 \\
 Wujie Shi\\
{\footnotesize\small School of Math.,  Suzhou
University, Suzhou 215006, Jiangsu, P.R. China}\\
{\footnotesize\small E-mail: wjshi@suda.edu.cn}\\
}

\begin{document}
\maketitle
\date{}

\begin{abstract}  The aim of  this paper is to classify the finite
nonsolvable  groups in which every irreducible character of even
degree vanishes on at most two conjugacy classes. As a corollary,
it is shown that $L_2(2^f)$ are the only nonsolvable groups in
which every irreducible character of even degree  vanishes on just
one conjugacy class.
\end{abstract}

\maketitle


\textbf{2000 MR Subject Classification}\,\,20C15

\section{Introduction}

For the simple groups $L_2(2^f)$ W. Feit given a characterization
using the centralizer of any involution in its Sylow 2-subgroups
(see \cite{F}). In this paper, we characterize $L_2(2^f)$ using
the number of even degree character zeros. For an irreducible
character $\chi$ of a finite group $G$, set $v(\chi)=\{g\in
G\,|\,\, \chi(g)=0\}$. Clearly, $v(\chi)$ is a union of some
conjugacy classes of $G$. An old theorem of Burnside asserts that
$v(\chi)$ is not empty for any nonlinear $\chi \in {\rm Irr}(G)$.
It makes sense to consider the structure of a finite group
provided that the number of character zeros in its character table
is small (see \cite{BK}, \cite{Q} for a few examples).  Y.
Berkovich and L. Kazarin \cite{BK} posed the following question:

\em Study the nonsolvable groups with the following ``property
(*):  every irreducible character of even degree has just one
zero, i.e., vanishes on just one conjugacy class''.  Is it true
that $L_2(2^f)$ are the only simple groups with this property? \em

Our observation is as follows.

\noindent {\bf Theorem A}\,\, \em  $L_2(2^f)$ ($f\geq 2$) are the
only nonsolvable groups with the property (*).\em

\smallskip In deed, we establish the following two theorems.

\noindent {\bf Theorem B}\,\, \em  Let $G$ be a finite group. Then
$G$ satisfies the property (*) if and only if $G$ is one of the
following groups: \em

(1)  \em  $G$ possesses a normal and Abelian Sylow $2$-subgroup;
\em

(2)  \em  $G$ is a Frobenius group with a complement of order $2$;
\em

(3)  \em  $G\cong SL(2,3)$; \em

(4)  \em  $G\cong L_2(2^f), f\geq 2$. \em

\bigskip

\noindent {\bf Theorem C}\,\, \em A finite nonsolvable group $G$
has the following ``property (**):\ every nonlinear $\chi\in {\rm
Irr}(G)$ of even degree has at most two zeros''  if and only if
$G\cong L_2(7)$ or $L_2(2^f)\times Z$ where $f\geq 2$ and
$|Z|=1,2$. \em

\bigskip

In this paper, $G$ always denotes  a finite group. For a subset
$A$ of $G$, let $k_G(A)$ be the minimal integer $l$ such that $A$
is a subset of a union of  $l$ conjugacy classes of $G$.  For
$N\triangleleft G$,  we put ${\rm Irr}(G|N)={\rm Irr}(G)-{\rm
Irr}(G/N)$; and for $\lambda\in {\rm Irr}(N)$, the inertia
subgroup of $\lambda$ in $G$ is denoted by $I_G(\lambda)$.

Let ${\rm Irr}_2(G)$ be the set of irreducible characters of $G$
with even degree.

We shall freely use the following facts: Let $N\triangleleft G$
and set $\overline{G}=G/N$.

{\rm (1)} For any $x\in G$, $\overline{x}^{\overline{G}}$ is a
union of some classes of $G$; furthermore,
$k_G(\overline{x}^{\overline{G}})=1$ if and only if $\chi(x)=0$
for any $\chi\in {\rm Irr}(G|N)$.

{\rm (2)} If $G$ has the property (*) or (**), then so has $G/N$.

\section{ Theorem C}

\begin{lem} \label{l201} Let $N<M$ be two normal subgroups of $G$ with
$k_G(M-N)=1$. Then $M$ is solvable.
\end{lem}
\pf Clearly, $M/N$ is an Abelian chief factor of $G$; also, all
elements in $M-N$ have the same order $p^d$, a power of some prime
$p$. For $x\in M-N$, since $M/N-N/N$ is also a $G/N$-class, we see
that $|C_{G/N}(xN)|=|C_G(x)|$. This implies that for every
$\chi\in {\rm Irr}(G|N)$, $\chi(x)=0$, and so that $p$ divides
$\chi(1)$ because $o(x)=p^d$ (see \cite[Lemma 1.1]{QSY}. By
\cite[Theorem D]{IK}$N$ is solvable, and so is $M$. \pfend

\begin{lem} \label{l202} For any nonlinear $\chi\in {\rm Irr}(G)$, we
have:

{\rm (1)}  If $G$ is nonsolvable and $k_G(v(\chi))\leq 2$, then
$\chi_{_{G'}}$ is irreducible, and

{\rm (2)} if $v(\chi)\subset N$ for some $N\triangleleft G$, then
$gcd(\chi(1), |G/N|)=1$.
\end{lem}

\pf  (1)\,\, Suppose that $\chi_{_{G'}}$ is reducible. Then
$\chi=\psi^G$, where $\psi\in {\rm Irr}(M)$ and $G'\leq M< G$.
This implies that  $\chi$ vanishes on $G-M$, and thus
$k_G(G-M)\leq 2$. By \cite[Theorem 2.2]{QSY} $G$ is solvable, a
contradiction.

(2)\,\, If $v(\chi)\subset N$, arguing as in the proof of
\cite[Lemma 3.14]{I} we conclude that $|\chi(x)|=1$ whenever $x\in
G-N$.  It follows from \cite[Lemma 1.1]{QSY} that $p
{\not|}\chi(1)$ for any prime divisor $p$ of $|G/N|$. Therefore
$gcd(\chi(1), |G/N|)=1$. \pfend

\begin{lem} \label{l203} Let $G$ be a non-Abelian simple group.
Then there exists $\chi\in {\rm Irr}_2(G)$ such that $\chi$ is of
$p$-defect zero for some prime divisor $p$ of $|G|$.\end{lem}

\pf It suffices to consider the non-Abelian simple group $G$ with
no irreducible character of $2$-defect zero.  By
\cite[Corollary]{W}, we may assume $G$ is isomorphic to one of the
following groups: $A_n$, $M_{12}$, $M_{22}$, $M_{24}$, $J_2$,
$HS$, $Suz$, $Ru$, $Co_1$, $Co_3$, or $B$.

If $G$ is isomorphic to $A_n, n\leq 8$ or one of the above
sporadic simple groups, we conclude the result by \cite{Atlas}.
Suppose that $G\cong A_n, n\geq 9$.  By \cite[Proposition]{MSW},
there is $\chi\in {\rm Irr}(G)$ such that $2p|\chi(1)$, where $p$
is the maximal prime not exceeding $n$. Clearly, $\chi$ is of
$p$-defect zero since $|G|_p=p$. \pfend

\begin{lem}\label{lnew} Let $G$ be a nonsolvable group with the
property (**). Then $G$ has the unique non-Abelian composite
factor.
\end{lem}

\pf By induction we may assume that $Sol(G)$, the maximal solvable
normal subgroup of $G$, is trivial. Let $N$ be a nonsolvable
minimal normal subgroup of $G$. If $N$ is not simple, then
$N=N_1\times \cdots \times N_s$ where each $N_i$ is a non-Abelian
simple group and $s\geq 2$. Let $\theta_i\in {\rm Irr}_2(N_i)$ be
of $p$-defect zero (Lemma \ref{l203}), and set
$\theta=\theta_1\times \cdots \times \theta_s$. Let $\chi_0$ be an
irreducible constituent of $\theta^G$,  let $x_1\in N_1, x_2\in
N_2$ be of order $p$ and $y_2\in N_2$ be of a prime order $q$
($q\not=p$). Clearly $\theta^g$ is of $p$-defect zero for any
$g\in G$, thus
$\theta^g(x_1)=\theta^g(x_1x_2)=\theta^g(x_1y_2)=0$. This implies
that $\chi_0(x_1)=\chi_0(x_1x_2)=\chi_0(x_1y_2)=0$. Since $x_1,
x_1x_2, x_1y_2$  lie in distinct conjugacy classes, we conclude a
contradiction. Suppose that $N$ is simple but $G/N$ is
nonsolvable. Since $Out(N)$ is solvable,  $C_G(N)$ is nonsolvable
and hence contains  a nonsolvable minimal normal subgroup $M$ as
$Sol(C_G(N))=1$. Set $T=M \times N$. Let $\psi\in {\rm Irr}_2(M)$
be $q$-defect zero and $\theta\in {\rm Irr}(N)$ be of $p$-defect
zero, where $q, p$ are prime divisors of $|M|$ and $N$
respectively. Let $x\in M$, $y\in N$ be of order $q,p$
respectively. Then for any irreducible constituent $\chi$ of
$(\psi \times \theta)^G$, we see that
$\chi(x)=\chi(y)=\chi(xy)=0$. Clearly, $x, y, xy$ lie in distinct
classes of $G$, a contradiction. \pfend

\begin{pro}\label{p201} Let $G$ be a nonsolvable group with the
property (**). If every  nontrivial quotient group of $G$ is
solvable, then $G\cong L_2(2^f), f\geq 2$ or $L_2(7)$.
\end{pro}

\pf Let $N$ be the unique minimal normal subgroup of $G$.  By
Lemma \ref{lnew} $N$ is a non-Abelian simple groups. Also $G/N$ is
solvable and $C_G(N)=1$.

\bigskip
Step 1. $N$ is a simple group of Lie type.

Suppose that $N\cong A_n$ for some $n\geq 8$.  Let $\pi$ be the
permutation character of $N$, and $\delta$ be the mapping of $N$
into $\{0,1,2,\cdots\}$ such that $\delta(g)$ is the number of
$2$-cycles in the standard  composition of $g$. Set
$$\lambda=\frac{(\pi-1)(\pi-2)}{2}-\delta, \,\,\,\,
\rho=\frac{\pi(\pi-3)}{2}+\delta.$$ \noindent By \cite[Ch.5,
Theorem 20.6]{H}, both $\lambda$ and $\rho$ are irreducible
characters of $N$. Clearly either $\lambda$ or $\rho$ is of even
degree. Let $\chi_0$ be an irreducible constituent of $\tau^G$,
where $\tau\in \{\lambda, \rho\}$ is of even degree.  By Lemma
\ref{l202}(1), $(\chi_0)_{_N}=\tau$.  For even $n$, set

$a_1=(1,\cdots , n-1)$, $a_2=(1,\cdots, n-2)(n-1, n)$,
$a_3=(1,\cdots, n-5)(n-4,n-3,n-2)$;

$b_1=(1,\cdots, n-3), b_2=(1,2,\cdots n-3)(n-2,n-1,n),
b_3=(1,\cdots, n-4)(n-3, n-2)$.

\noindent For odd $n$, set

$a_1=(1,\cdots, n-2)$, $a_2=(1,\cdots,n-4)(n-3,n-2,n-1)$,
$a_3=(1,\cdots, n-5)(n-4,n-3)$;

$b_1=(1,\cdots,n)$, $b_2=(1,\cdots,n-3)(n-2,n-1)$, $b_3=(1,\cdots,
n-6)(n-5, n-4,n-3)$.

\noindent We see that $\lambda(a_i)=0=\rho(b_i)$ for any
$i=1,2,3$.  Therefore, either
$\chi_0(a_1)=\chi_0(a_2)=\chi_0(a_3)=0$ or
$\chi_0(b_1)=\chi_0(b_2)=\chi_0(b_3)=0$. Observe that $a_1, a_2,
a_3$ (or $b_1, b_2, b_3$) lie in distinct classes of $G$. We
conclude a contradiction.

Suppose that $N\cong A_7$ or  one of the sporadic simple groups.
For the case when  $G>N$, then $|G/N|=2$. It implies  by Lemma
\ref{l202} that for any non-principal $\theta\in {\rm Irr}(N)$,
$\theta$ is extendible to $\chi\in {\rm Irr}(G)$,  and that
$k_G(v(\chi)\cap N)=1$ whenever $\theta(1)$ is even. By
\cite{Atlas}, we conclude a contradiction. For the case when
$G=N$, we also conclude a contradiction by \cite{Atlas}.

Note that $A_5\cong L_2(4)\cong L_2(5)$, $A_6\cong L_2(9)$. By the
classification theorem, $N$ must be a simple group of Lie type.

\bigskip

Step 2. If  $G=L_2(q)$ for some odd $q=p^f>5$, then $G\cong
L_2(7)$.

Since  $N$ is one of the simple groups of Lie type, by \cite{W}
$N$ has an irreducible character $\chi_0$ of $2$-defect zero. Let
$\sigma$ be an irreducible constituent of $\chi_0^G$. Observe that
$\chi_0^g(x)=0$ for any $g\in G$ and any $x\in N$ of even order.
It follows that $\sigma(x)=0$ whenever $x\in N$ is of even order.
Let $P\in Syl_2(N)$, and $\Delta=\cup_{g\in G}(P^g-\{1\})$.

We claim that if  $G=L_2(q)$ for some odd $q=p^f>5$, then $G\cong
L_2(7)$. If $\sigma(1)\not=p^f+1$ where $\sigma \in {\rm Irr}(G)$
is of $2$-defect zero, then set $\eta\in {\rm Irr}(G)$ with degree
$p^f+1$. Let $C$ be a single cycle of $G$, and $\Xi=\cup_{g\in
G}(C^g-1)$. Clearly, for any $v\in \Xi$, either $\sigma(v)$ or
$\eta(v)=0$. This implies that $k_G(\Xi)\leq 4$. Since
$k_G(\Xi)=(q-1)/4$ (see \cite[II, Theorem 8.5]{H}), we have that
$q=7,9,11,13$. By \cite{Atlas}, we conclude that $q=7$ and $G\cong
L_2(7)$.

\bigskip
Step 3. If $P$ is non-Abelian, then $G\cong L_2(7)$.

In this case, we have $v(\sigma)=\Delta \subset N$ and
$k_G(v(\sigma))=2$. By Lemma \ref{l202}, $|G/N|$ is odd and
$\sigma_{_N}=\chi_0$. Therefore $\sigma$ is of $2$-defect zero,
and $\sigma(x)=0$ for any $x\in G$ of even order. This implies
that $P\in Syl_2(G)$ and $C_G(t)$ is a $2$-group for any
involution $t$. By \cite[III, Theorem 5]{S} and $P$ is
non-Abelian, $G$ is one of the following groups: $Sz(q),
q=2^{2m+1}$, $L_2(q)$ where $q$ is a Fermat prime or Mersenne
prime, $L_3(4)$, $L_2(9)$ or $M_{10}$.

By \cite{Atlas} $M_{10}$ and $L_3(4)$ do not have the property
(**). Note that all elements of order $4$ in $Sz(2^{2m+1})$ forms
two conjugacy classes , which can be easily verified by
\cite[Ch.XI, Theorem 3.10]{HB}. Therefore $G\cong Sz(2^{2m+1})$ is
not the case. Now by step 2, we conclude that $G\cong L_2(7)$.

\bigskip

Step 4. If $P$ is Abelian, then  $G\cong L_2(2^f), f\geq 2$.

Since $P$ is Abelain,  by \cite[XI, Theorem 13.7]{HB}, $N$ is one
of the following groups: $L_2(2^f)$, $L_2(q)$ where $q=3,5(mod\,\,
8)$ or $^2G_2(q)$ where $q=3^{2m+1}$. Recall that $\sigma(x)=0$
whenever $x\in N$ is of even order.

Suppose that  $N\cong \,\,^2G_2(q)$.  Then all elements of even
order in $N$ lie in at least  three  classes of $G$ (see
\cite[Theorem 13.4]{HB}), a contradiction.

Therefore $N\cong L_2(q)$. Then $Aut(N)=N\la \phi, \delta\ra$,
where $\la \phi\ra$ is the group of field automorphisms of $N$,
$\la \delta\ra$ is the group of diagonal automorphisms of $N$.
Note that $\phi$ and $\delta$ commute modulo $Inn(N)$, and that if
$q$ is even then $\la\delta \ra$ is trivial.

Suppose that $N=L_2(q)$ where $q>5, q=3,5(mod\,\,8)$. Let $I=N\la
\delta\ra \cap G$ and $\theta\in {\rm Irr}(N)$ of degree $q-1$.
Arguing as the proof of \cite[Theorem 2.7]{LW}, we conclude that
$\theta$ is extendible to $\psi\in {\rm Irr}(I)$ and that $\psi^G$
is irreducible.  It follows by Lemma \ref{l202} that $G=I$. In
particular, $G=N$ or $G=N\la \delta\ra$. For the case when $G=N$,
step 2 already yields a contradiction. Suppose that  $G=N \la
\delta \ra \cong PGL(2,q)$. Let $\theta, \theta_1\in {\rm Irr}(N)$
with degree $q-1$, $q+1$ respectively. By Lemma \ref{l202},
$\theta, \theta_1$ are extendible to $\mu$ and $\mu_1\in {\rm
Irr}(G)$ respectively, and $k_G(v(\mu)\cap N)=k_G(v(\mu_1)\cap
N)=1$. Observe that either $(q-1)/2$ or $(q+1)/2$ is equal to $2k$
for some $k>1$, and that $\mu(a)=\mu_1(b)=0$ whenever $a, b\in N$
with $2o(a)|(q-1), 2o(b)|(q+1)$. It follows that either
$k_G(v(\mu)\cap N)$ or $k_G(v(\mu_1)\cap N)$ is greater than $1$,
a contradiction.

Therefore $N\cong L_2(2^f)$. Suppose that $G>N$. Then  $G=G\cap
N\la \phi\ra$. For $\theta\in {\rm Irr}(N)$ with degree $2^f-1$,
$\theta$ is $\phi$-invariant and hence induces to an irreducible
character of $G$. This implies by Lemma \ref{l202} that  $G/N$ is
a cyclic group of odd order. Recall that $\chi_0\in {\rm Irr}(N)$
is of degree $2^f$, and $\sigma$ is an extension of $\chi_0$ to
$G$. Note that since $\Delta=\cup_{g\in G}(P^g-1)$ is a class of
$N$, it forces  $\Delta$ to be also a class of $G$. This implies
that $|C_G(t)|=|G/N||P|$ for  $t\in P-\{1\}$, and so that
$C_G(t)=PA$, where $A\cap N=1$, $A\cong G/N$. Observe that
$\sigma(g)=0$ whenever $g\in G$ is of even order. It follows that
$|G/N|$ is a odd prime $q$ and that $\Theta$, the set of elements
of order $2q$, forms a class of $G$. Let $w\in A$ be an element of
order $q$, and $y=wt$. Since $\Theta$ is a class of $G$, all
subgroup of order $2q$ are conjugate. Therefore,
$|\Theta|=|G:N_G(\la y\ra)|(q-1)$. It is easy to see that
$C_G(y)\leq N_G(\la y\ra)\leq C_G(t)\cap C_G(w)=C_G(y)$.
  Then
$$|G:C_G(y)|(q-1)=|G:N_G(\la y\ra)|(q-1)=|\Theta|=|G: C_G(y)|, $$
\noindent a contradiction. \pfend

\begin{lem} \label{l204} Let $N\triangleleft G$ and $H/N$ be a Hall $\pi$-subgroup of
$G/N$. If  $\eta\in {\rm Irr}(H)$ induces to an irreducible
character $\chi$ of $G$,  then $\chi(x)=0$ for any $\pi'$-element
$x\in G-N$.
\end{lem}

\pf  It follows directly from the definition of  induced
character. \pfend

\begin{lem} \label{l205} Suppose that $H$ is a subgroup of $G\cong L_2(2^f), f\geq
3$. Let $P_1\in Syl_2(H)$, $P\in Syl_2(G)$ be such that $P_1\leq
P$. Then, we have

{\rm (1)} If $|P:P_1|=2$, then $H=P_1$.

{\rm (2)} If $P=P_1$, then $P\leq H\leq N_G(P)$.
\end{lem}

\pf It is enough to investigate the maximal subgroups of
$L_2(2^f)$ (see \cite[II, Theorem 8.27]{H}). \pfend

\bigskip
\noindent{\bf \em Proof of Theorem C\em}\,\, It suffices to prove
that if $G$ is a nonsolvable group with the property (**), then
$G\cong L_2(7)$ or $ L_2(2^f)\times Z$ where $f\geq 2$ and
$|Z|\leq 2$.

Let $N\triangleleft G$ be maximal such that $G/N$ is nonsolvable.
By  Lemma \ref{lnew} and Proposition \ref{p201}, $N$ is solvable
and  $G/N\cong L_2(2^f)(f\geq 2)$ or $L_2(7)$. Set
$\overline{G}=G/N$, and $\sigma\in {\rm Irr}(\overline{G})$ be of
degree $2^f$. Let $P>N$ be such that $\overline{P}\in
Syl_2(\overline{G})$, and let $\Delta:=\cup_{g\in G}(P^g-N)$.
Note that for any $f\geq 2$, if $L_2(2^f)$ has nontrivial Schur
multiplier, then $f=2$.

We first claim that if $\overline{G}=L_2(7)$ then $N=1$. If else,
to see a contradiction we may assume that $N$ is minimal normal.
Observe that $\Delta\subseteq v(\sigma)$ and
$k_{\overline{G}}(\overline{\Delta})=2$. It follows  that
$k_G(\Delta)=2$, and  that $C_G(t)$ is a $2$-group,
$|C_G(t)|=|C_{\overline{G}}(\overline{t})|$ for any $t\in \Delta$.
In particular,  there is $t_0\in \Delta$ such that $|C_G(t_0)|=4$.
This implies that a Sylow $2$-subgroup of $G$ possesses a cyclic
subgroup of index $2$ (see \cite[Lemma 1.3]{QSY}). Thus a Sylow
$2$-subgroup of $N$ is cyclic, and so either $|N|=2$ or $|N|=q^r$
for some odd prime $q$. If $|N|=2$, then $G\cong L_2(7)\times
Z(G)$, which is clear not the case.  If $|N|=q^r$ is odd, then $P$
is a Frobenius group with a complement  $\cong P/N\cong D_8$,
which is also impossible.

We claim that if $\overline{G}=L_2(2^f)$ and $G'<G$, then $G'\cap
N=1$ and $G\cong L_2(2^f) \times Z$ where $|Z|=2$. It is easy to
see that $G/G'\cap N\cong L_2(2^f) \times Z$ where $|Z|=2$. If
$K:=N\cap G'>1$, to see a contradiction we may assume $K$ is
minimal normal. Note that $P/K=P_1/K \times Z/K$ is elementary
Abelian and of order $2^{f+1}$. Since $2=k_{G/K}(\Delta)\leq
k_G(\Delta)\leq k_G(v(\sigma))\leq 2$, we conclude that for any
$t\in P_1-K \subset P-N$, $|C_G(t)|=|C_{G/K}(tK)|=2^{f+1}$. Now
for any $\mu\in {\rm Irr}(G|K)$, since $\mu$ vanishes on $\Delta$,
we see that $2|\mu(1)$ (see \cite[Lemma 1.1]{QSY}), and so
$$v(\mu)=\Delta.$$ \noindent If $N$ is of odd order, then $P_1$
acts Frobeniusly on $N$, which is clearly impossible. Suppose that
$N$ is an $2$-group. Let non-principal $\lambda\in {\rm Irr}(K)$
be such that $\lambda$ is $P$-invariant, and $\mu$ be some
irreducible constituent of $\lambda^G$. If $I_G(\lambda)=G$, then
$K\leq Z(G)$, and so $G=L_2(2^f) \times K$ or $G\cong
SL(2,5)\times Z_1$ ($|Z_1|=2$), which is impossible. If
$I_G(\lambda)<G$, then $P\leq I_G(\lambda)\leq N_G(P)$, and then
$\mu(x)=0$ whenever $o(x)|2^f+1$ (see Lemma \ref{l204}), which
contradicts the claim: $v(\mu)=\Delta$.

In what follows, we need only to show that $N=1$ provided that
$\overline{G}\cong L_2(2^f)$ and $G=G'$. If else, to see a
contradiction we may assume that $N$ is minimal normal, and so is
an elementary Abelian  $q$-group for some prime $q$.

For any non-principal $\lambda\in {\rm Irr}(N)$, if $\lambda$ is
$G$-invariant, then $N=Z(G)$, and since $G=G'$ we conclude that
$N$ is a subgroup of the Schur multiplier of $L_2(2^f)$, and so
that $G\cong SL_2(5)$. By \cite{Atlas}, $SL(2,5)$ does not satisfy
the property (**), a contradiction. Therefore, $I_G(\lambda)<G$
for any non-principal $\lambda\in {\rm Irr}(N)$. Using the same
argument as in above, we also conclude that $|N|>2$.

Let $P<H<G$ be such that
$\overline{H}=N_{\overline{G}}(\overline{P})$. We claim that
\em``there is $\lambda_0\in {\rm Irr}(N)$ such that
 $I_G(\lambda_0)\leq H$, and $\lambda_0^G$ has an irreducible
constituent $\chi_0$ of even degree''.\em

Case a. $q=2$, i.e., $N$ is a $2$-group.

Since $\Delta\subseteq v(\sigma)$, $k_G(\Delta)=1,2$. Note that
$C_G(x)$ is a $2$-group and that $|C_G(x)|\leq 2^{f+1}$ since
$k_G(\Delta)\leq 2$. This implies that for any $x\in P-N$,
$$2^f=|P/N|\leq |P/P'|\leq |P/P'|+\sum_{\eta\in {\rm Irr}(N),
\eta(1)>1}|\eta(x)|^2=|C_P(x)|\leq |C_G(x)|\leq 2^{f+1},$$
\noindent and so $|C_P(x)|=2^f$ or $2^{f+1}$, $|N : P'|=1,2$.

Suppose that $N=P'$ and let $P$-invariant non-principal
$\lambda_0\in {\rm Irr}(N)$. By Lemma \ref{l205}, we  see that
$P\leq I_G(\lambda_0)\leq H$, also that $\lambda_0^G$ has an
irreducible constituent of even degree.

Suppose that  $|N: P'|=2$.  By the above inequality, we conclude
that any nonlinear irreducible character of $P$ must vanish on
$P-N$. Note that since $|N|>2$, $P'>1$. Let $P'/E_1$ be a
principal factor of $P$ and let $ E(\geq E_1)\triangleleft P$ be
maximal such that $P/E$ is non-Abelian. Let $V/E=Z(N/E)$ and let
$\eta_0$ be a nonlinear irreducible character of $P/E$. Since $P-N
\subseteq v(\eta_0)=P-V$,  $V\leq N$. Observe that $V/E$ is cyclic
by \cite[Lemma 12.3]{I} but $V\leq N$ is elementary Abelian, it
follows that $V/E$ has order $2$,  and so that $V=N$ or $P'$. If
$V=N$, then $(\eta_0)_N=\eta_0(1)\lambda_0$. Thus $\lambda_0$ is
$P$-invariant, and this $\lambda_0$ works for our claim.  If
$V=P'$, then $(\eta_0)_{P'}=2^d\gamma$, where $2^{2d}=2^{f+1}$.
Since $\eta_0$ vanishes on $P-P'$, it is easy to check that
$\eta_N=2^{d-1}(\lambda_0+\lambda_1)$, where $\lambda_0$,
$\lambda_1$ are distinct linear characters of $N$. Therefore $|P:
I_P(\lambda_0)|=2$. Note that $|P/N|\geq 8$ since $f+1=2d$. Now by
Lemma \ref{l205}, we see that such $\lambda_0$ also works for the
claim.

Case b. $q>2$, i.e., $N$ is of odd order.

In this case, since $P/N$ is elementary Abelian, it is easy to
conclude that there is non-principal $\lambda_0\in {\rm Irr}(N)$
such that $|P: I_P(\lambda_0)|=2$. Thus, if  $2^f\geq 8$, then
such $\lambda_0$ works for our claim by Lemma \ref{l205}.

For the case when $2^f=4$, we also conclude that $2$ divides $|
G:I_G(\lambda_0)|$ since $\overline{P}$ is a T.I subgroup of
$\overline{G}$.  Now if the claim  fails, then $K/N :=
I_G(\lambda_0)/N \cong D_{10}$ or $S_3$.  Set
$\lambda_0^K=e_1\theta_1+\cdots e_s\theta_s$, where $\theta_1,
\cdots, \theta_s\in {\rm Irr}(K)$ are distinct. Observe that
$e_i=\theta_i(1)$ divides $|K/N|$, that $|K/N|=e_1^2 +\cdots
+e_s^2$,  and that if $1\in \{e_1, \cdots, e_s\}$ then $2\in
\{e_1, \cdots, e_s\}$ by \cite[Corrollary 6.17]{I}. It follows
that $2\in \{e_1, \cdots, e_s\}$. Therefore, there is $\omega \in
{\rm Irr}(O^2(K))$ such that $\omega^G$ is irreducible. Now by
Lemma \ref{l204}, it is easy to see that $k_G(v(\omega^G))\geq 3$,
a contradiction. This complete the proof of the claim.

Let $\lambda_0, \chi_0$ be as in the claim. Then $\chi_0=\omega^G$
for some $\omega\in {\rm Irr}(H)$. Let $\Xi$ be the set of
elements outside $N$ of order divisible by $2^f+1$. By Lemma
\ref{l204}, $\Xi \subseteq v(\chi)$ and hence $k_G(\Xi)\leq 2$.
Note that $k_G(\Xi)\geq k_{\overline{G}}(\overline{\Xi})=2^{f-1}$.
This implies that $f=2$, $\overline{G}=L_2(4)$.

Let us investigate $I_G(\lambda_0)\leq H$. We see that either
$I_G(\lambda_0)/N$ is a $2$-group, $3$-group, or
$I_G(\lambda_0)=H\cong A_4$.

Suppose that  $I_G(\lambda_0)/N$ is a $2$-group or a $3$-group. By
Lemma \ref{l204}, we can easily conclude a  contradiction.

Suppose that  $I_G(\lambda_0)=H$. Since $\lambda_0^G$ has an
irreducible constituent $\chi_0$ of even degree, we see that
$\lambda_0^P$ has an irreducible constituent of even degree. This
implies that $\lambda_0^P=2\psi$ where $\psi\in {\rm Irr}(P)$ is
of degree $2$. Let $\psi_1$ be any irreducible constituent of
$(\chi_0)_P$ and let $\lambda_1$ be an irreducible constituent of
$(\psi_1)_{_N}$. Clearly, $\lambda_1=\lambda_0^g$ for some $g\in
G$. Note that $I_G(\lambda_1)=H^g$, and that either $P\leq H^g$ or
$P\cap H^g=N$ since $\overline{P}$ is a T.I set of $\overline{G}$.
It follows that either $\lambda_1=\lambda$ or $I_P(\lambda_1)=N$.
Therefore, $\psi_1$ always vanishes on $P-N$. This implies that
$\chi_0$ vanishes on $\Delta\cup \Xi$, which contradicts the
property (**). This completes the proof of Theorem C.

\section{Theorem A and Theorem B}

{\bf \em Proof of Theorem A\em}\,\, It is a direct consequence of
Theorem C.

\bigskip

\noindent{\bf \em Proof of Theorem B\em}\,\, It is obvious that
all the groups listed in Theorem B have the property (*). Note
that if ${\rm Irr}_2(G)$ is empty, then $G$ possesses an Abelian
and normal Sylow $2$-subgroup.  Now we need only to show that if
$G$ is a solvable group with the property (*) and ${\rm Irr}_2(G)
\not=\emptyset$, then either $G$ is a Frobenius group with a
complement of order $2$, or $G\cong SL(2,3)$.

Suppose first that there is some $\chi\in {\rm Irr}_2(G)$ such
that $\chi_{_{G'}}$ is reducible. Then $\chi=\lambda^G$ for some
$\lambda\in {\rm Irr}(H)$, where $G'\leq H<G$. This implies that
$\chi$ vanishes on $G-H$, and so $k_G(G-H)=1$. It is easy to
conclude in this case that $G$ is a Frobenius group with a
complement of order $2$.

Now it suffices to prove that $G\cong SL(2,3)$ provided that
$\chi_{_{G'}}$ is irreducible for any $\chi\in {\rm Irr}_2(G)$.
Observe that for any $\chi\in {\rm Irr}_2(G)$,  $v(\chi)\cap G'$
can not be empty, and it follows by Lemma \ref{l202} that
$gcd(\chi(1), |G/G'|)=1$. In particular, $|G/G'|$ is odd.

Let $E\triangleleft G$  maximal be such that $G/E$ is non-Abelian.
By \cite[Lemma 12.3]{I} $G/E$ is a $p$-group or a Frobenius group.
If $G/E$ is a $p$-group and let nonlinear $\psi\in {\rm
Irr}(G/E)$, since $gcd(\chi(1), p)=1$ for any $\chi\in {\rm
Irr}_2(G)$, we get that $\chi\psi\in {\rm Irr}_2(G)$, which is
impossible because $gcd((\chi\psi)(1), p)>1$. Therefore $G/E$ is a
Frobenius group with a kernel $N/E$ and a cyclic complement, and
thus $gcd(\chi(1), |G/N|)=1$ whenever $\chi\in {\rm Irr}_2(G)$.
Clearly, for any $\tau\in {\rm Irr}_2(N)$, $\tau$ is extendible to
some $\chi\in {\rm Irr}_2(G)$, and thus \cite[Theorem 12.4]{I}
implies that both $\chi$ and $\tau$ vanish on $N-E$, so
$k_G(N-E)=1$. Observe that $k_{G/E}(N/E-E/E)=1$. It follows that
$|N/E|=1+|G/N|$. Since $|G/N|$ is odd, we see that $N/E$ is a
$2$-group.  Set $|N/E|=2^r$ and we have
$2^r=|C_G(t)|=|C_N(t)|=|N/E|+\sum_{\eta\in {\rm
Irr}(N|E)}|\eta(t)|^2$ for any $t\in N-E$. This implies that: \em
``for any $\eta\in {\rm Irr}(N|E)$, $\eta$ must vanishes on $N-E$
,  so $2|\eta(1)$ (\cite[Lemma 1.1]{QSY}), and hence $\eta$ is
extendible to $G$''. \em

Clearly we may assume that $E>1$ and set $E/F$ be a principal
factor of $G$. If $E/F$ is not a $2$-group, then the above fact
implies that $N/F$ is a Frobenius group with the kernel $E/F$, and
then $N/E$ is of order $2$ which is impossible. Thus $E/F$ is a
$2$-group. Let us investigate the quotient group $G/F$. Let
$T\cong G/N$ be a Hall $2'$-subgroup of $G/F$.  We see that  $K$
acts non-trivially on $N/F$ and fixes every nonlinear irreducible
character of $N/F$. By \cite[Lemma 19.2]{MW}, we conclude that
$N'/F=E/F\leq Z(G/F)$. Now \cite[Ch.5, Theorem 6.5]{G} implies
that $2^r-1=|G/N|$ divides $2^e+1$ for some integer $e\leq r/2$.
This yields that $2^r=4$, and so $G/F\cong SL(2,3)$.

To finish the proof of Theorem B, it suffices  to show that $F=1$.
If else, towards a contradiction we may assume that $F$ is a
minimal normal subgroup of $G$. Suppose that  $F$ is a $2$-group.
Since $|C_G(t)|=4$ for any $t\in N-E$, there is $x\in N-E$  of
order $|N|/2\geq 8 $ (\cite[Lemma 1.3]{QSY}), which is clearly
impossible. Suppose that $F$ is a $q$-group for some odd prime $q$
and set $P\in Syl_2(N)$. Since $|C_G(t)|=4$ for any $t\in N-E$, we
see that $C_P(x)\leq P'$ for any $1\not=x\in F$. It follows by
\cite[Lemma 19.1]{MW} that $P$ is either cyclic or isomorphic to
$Q_8$, which is also impossible. The proof of Theorem B is
complete.

\smallskip

By Theorem C, it is not difficult to deduce  the following known
result (\cite{BGC}).

\noindent{\bf Proposition} \em Let $G$ be a finite non-Abelian
group in which every irreducible character has at most two classes
of zeros. Then $G$ is one of the following groups:\em

(1) \em$G$ is Frobenius group with a complement of order $2$ or
$3$;\em

(2) \em $G/Z(G)$ is a Frobenius group with a complement of order
$2$, and $|Z(G)|=2$;\em

(3) \em $G\cong S_4$;\em

(4) \em$G\cong L_2(4)$ or $L_2(7)$.\em

\nopagebreak

\end{document}